\documentclass{article} 
\usepackage{amsmath}
\usepackage{paralist}
\usepackage{subcaption}
\usepackage{lscape}
\usepackage[misc]{ifsym}
\usepackage{epsfig} 
\usepackage{epstopdf} 
\usepackage[colorlinks=true]{hyperref}
\hypersetup{urlcolor=blue, citecolor=red}
\usepackage{hyperref}
\allowdisplaybreaks

\newtheorem{theorem}{Theorem}[section]

\newtheorem{remark}[theorem]{Remark}

\title{Management Strategies for Hydropower Plants \\ - a simple deterministic dynamic programming approach} 
\author{Marcus Olofsson\thanks{Corresponding author. e-mail: bakkenolofsson@gmail.com}   \\ \small{Umeå University, Sweden} \\ \small{Odensalaskolan Östersund, Sweden} \and André Berg \\ \small{Umeå University, Sweden}  \and Niklas L.P. Lundström \\ \small{Umeå University, Sweden} } 
\date{}
\begin{document}

\maketitle
\begin{abstract}
We use a dynamic programming approach to construct management strategies for a hydropower plant with a dam and a continuously adjustable unit. Along the way, we estimate unknown variables via simple models using historical data and forecasts. Our suggested scheme achieves on average 97.5 $\%$ of the theoretical maximum (optimal strategy when knowing the future) with small computational complexity. We also apply our scheme to a Run-of-River hydropower plant and compare the strategies and results to the much more involved PDE-based optimal switching method studied in \cite{OOL22}; this comparison shows that our simple approach may be preferable when the 
forecast is good enough.
\end{abstract}


\section{Introduction}
Even though operation of hydropower plants has been frequently studied during decades \cite{singh, Queiroz}, 
there is still a need for improving management of hydropower. As a green alternative for balancing volatile energy sources, small hydropower plants are increasingly important for a sustainable and stable electricity grid and effective management strategies are therefore necessary to increase their economic appeal. As an example, \cite{C_etal15} conclude that small hydropower is one of the most important impetuses for the development of China's power industry, and that the dispatching of their many small hydropower plants is lagging behind the development of the power grid and other power sources.

Managing hydropower is however a nontrivial task, even in the case of only a single power plant with one reservoir. When running a unit (a turbine connected to a generator), 
the reservoir is naturally drained of water and, if the inflow is insufficient, the head is lowered leading to less electricity generated per $m^3$ of water used.
Intuitively, one would therefore argue to keep the reservoir as close to full as possible, while minimizing the risk that it overflows (to avoid losing water without producing any electricity). However, as the cost of moving between production modes is often non-negligible, it might actually be optimal to allow some spillage of water to avoid paying this cost. 
Non-negligible costs for switching between production modes come from the fact that starting and stopping units
induces wear and tear on the machines and may also require intervention from
personnel \cite{NS97}. Each start and stop also involves a small risk, e.g., 
the major breakdown in the Akkats hydropower plant
(Lule river, Sweden) 2002 was caused by a unit being stopped too quickly,
resulting in rushing water destroying the foundation of both the turbine and the
generator \cite{Y10,Y18}.

It is a classical technique in optimization to use dynamic programming and backward induction to determine optimal decision sequences, so also in hydropower management.
The key idea is to answer the question: ''What is the best decision at this point, assuming that all my future actions will be optimal?" This method is extremely useful when all parameters are deterministic or to find the optimal decision \textit{in hindsight}, when the outcome of any random event is determined. In reality, however, systems typically involve randomness and there is not enough information available to determine the best decision without knowledge of the future. 

In this paper we overcome this lack of information by replacing the unknown random variables with estimates based on appropriate models, historical data, and forecasts, and by applying a dynamic programming technique to produce management strategies for a hydropower plant with a dam and a continuously adjustable unit. We use our estimates to construct \textit{approximate} optimal strategies and base our decisions on these approximations. When paired with short-term forecasts, this turns out to be a very efficient way that achieves close to optimal strategies with small modelling and computational effort. This method eludes deep mathematical theory that might be out of reach to practitioners and completely circumvents the need for simulations and numerical solutions of differential equations as used by the authors in \cite{OOL22}. 

The huge existing body of scientific literature on optimization methods to reservoir operation problems makes a state-of-the-art review difficult here, but we mention some contributions with no ambition or claim to being exhaustive. 
Three case studies where considered in \cite{DynProg86} for instantaneous, hourly, and monthly time frames. The study considers two linked power plants as well as cost of imported power supply.
State incremental dynamic programming was developed for multi-reservoir systems in \cite{YSA05}, in which a random file
access method was used for reaching data to cope with the curse of dimensionality.
In \cite{IDynProg14}, 
a concave approximation was derived and tested
on case studies of long-term hydropower scheduling, 
showing that the computation time increases linearly in
accordance with the number of storage intervals, whereas standard dynamic programming shows a quadratic increase.
In  \cite{IDynProg17}, a hybrid min-max dynamic programming model was formulated for peak operation of hydropower systems. 
The approach was tested on a large-scale
hydropower system in China, indicating satisfactory performance in reducing peak loads in the system.
In \cite{DFHD21} the flexibility that could be provided by large hydropower reservoirs in West Africa, to cope with planned
future solar and wind energy generation in the region, was investigated. 
%
%
Concerning review papers, an argument-driven classification and comparison of reservoir operation optimization methods can be found in \cite{DWP19}.
Heuristic programming methods, evolutionary and
genetic algorithms, as well as multi-objective optimization are discussed in \cite{rew04}, along with application of neural networks and fuzzy rule-based systems for inferring reservoir system operating rules.
A review on how the operation design problem is formulated, rather than solved, can be found in \cite{rew21}, including classification of over 300 studies published over the last years into distinctive categories depending on the adopted problem formulation.
Moreover, there is an extensive literature on how to improve and rationalize dynamic programming algorithms for managing hydropower, see, e.g., \cite{IDynProg14, IDynProg17} and the references therein, as well as other optimization techniques and models \cite{bok, rew04, DWP19, rew21, OOL22, YT22, highdim}. 

We contribute to the literature by showing how rudimentary optimization techniques together with simple mathematical models for river flow can yield very good results; the suggested production schemes performs remarkably well, averaging $97.5$~\% of the theoretical maximum (optimal strategy when knowing the future) over the years 2015-2022 when studying management of an example (fictitious) hydropower plant using real flow data from the northern parts of Sweden.
A possible explanation for this high performance is our suitable way of merging historical flow with short-term forecast in the simple flow model presented in Section \ref{sec:river-model}.

The rest of the paper is organized as follows.  
Section \ref{sec:method} outlines our suggested optimization scheme, including the dynamic programming approach,
river flow model, and the modeling of our power plant. 
The results from testing our scheme is presented in Section \ref{sec:results}, 
and in Section \ref{sec:oldsetup} we apply our scheme to a Run-of-River hydropower plant to compare the method with the more involved PDE-based method of \cite{OOL22}. 
We end in Section \ref{sec:discussionanconclusions} with a discussion of our findings and suggestions for further research.

\section{Problem setup and method}
\label{sec:method}

The objective in our problem is to manage a production facility to maximize its profit. More precisely, the manager must continuously choose between different modes of production, each with different profitability depending on some random process $X_t$. However, each change in production induces a cost and these costs are deducted from the total profit. This means that the optimal strategy is \textbf{not} to always switch to the mode with momentarily highest payoff. 

In the context of hydropower, this amounts to maximizing the profit from the electricity generated over a specific period $[t,T]$.  
More precisely, we want to maximize
\begin{equation} \label{eq:tooptimize}
\int_t^T \phi_{\mu_s}(s,Q_s,H_s,P_s)  \, ds - \sum_{t \leq \tau_i \leq T} c_{\xi_{i-1},\xi_i}(\tau_i,Q_{\tau_i},H_{\tau_i},P_{\tau_i}),
\end{equation}
by finding an optimal sequence of (random) time points $\tau_i$ at which we move from production mode $\xi_{i-1}$ to $\xi_i$. It is convenient to associate to each such sequence a (random) function $\mu_s$, indicating the current mode of production at time $s$, and we will move between these notations throughout the text without further notice. (In fact, we use both notations already in \eqref{eq:tooptimize}.) In the above display,
\begin{equation*}
\phi_{\mu_s}(s,Q_s,H_s,P_s)   
\end{equation*}
is the running payoff of the plant at time $s$ when in mode $\mu_s$, with water flow $Q_s$, reservoir head $H_s$, and electricity spot price $P_s$. Going forward, the variables $(q,h,p)$ indicate the current value of these stochastic processes, i.e.,  $Q_t=q,H_t=h,$ and $P_t=p$. The cost of moving from production mode $i$ to production mode $j$ is denoted $c_{ij}$. These costs occur due to, e.g., wear and tear of the components or the risk of failure when changing production mode.

\subsection{Dynamic programming for hydropower production}
\label{sec:method3}
The strategies constructed in this paper are based on dynamic programming paired with historical estimation and forecasts of water flows.
To put this method of optimization into the current setting, label the different modes of production $i \in \{1,2,\dots,m\}$. Let $V_i(t,x)$ be the optimal profit at time $t$ given that we are currently in production mode $i$ and that $X^{x,\alpha}_t:=(Q^q_t, H^{h,\alpha}_t, P^p_t)=(q,h,p)=:x$. The superindex $\alpha$ is used to stress that the process $X^{x,\alpha}_t$ is controlled in the sense that the reservoir level $H^{h,\alpha}$ depends on the amount of water used for production. {We will typically drop this superindex in favour of an easier notation. When in mode $i$ at time $t$, the total payoff from staying in mode $i$ until $t+\Delta t$ is
$$\phi_i(t,x) \cdot \Delta t + V_i(t+\Delta t, X^{x}_{t+\Delta t})$$
whereas switching to mode $j$ gives total payoff
\begin{equation} \label{eq:optimal}
\phi_j(t,x) \cdot \Delta t + V_j(t+\Delta t, X^x_{t+ \Delta t}) - c_{ij}(t,x).
\end{equation}
Therefore, the optimal decision is to choose whatever action that maximizes this output, i.e, to maximize \eqref{eq:optimal} over $j\in \{1,2,\dots,m\}$. With the value of acting optimally in the future $V_i(t+\Delta t,x)$, $i \in \{1,2,\dots,m\}$,  given, the optimal value $V_i(t,x)$ must therefore satisfy 
\begin{equation} \label{eq:optimalvalue}
V_i(t, x) = \max_{j \in \{1,2,\dots,m\}} \left\{\phi_j(t,x) \cdot \, \Delta t+ V_j(t+\Delta t, X^x_{t+\Delta}) - c_{ij}(t,x)\right\}.  
\end{equation}
If the terminal value $V_i(T,x)$ is known, we can thus work recursively backwards to find $V_i(t,x)$ for all $(t,x)$ and with $\{V_1, \dots, V_m\}$ known, the optimal decision $j^\ast$ in \eqref{eq:optimal} is given by
\begin{equation}\label{eq:strategy}
j^\ast = \arg \max_{\hspace{-0.65cm}	j \in \{1,2,\dots,m\}}  \left\{ \phi_j(t,x) \cdot \Delta t + V_j(t+\Delta t, X^x_{t+\Delta t}) - c_{ij}(t, x)\right\},
\end{equation}
where we assume $c_{ii}\equiv0$.

In applications, the assumption that $X_t$ is deterministic typically fails and the value of $X_{t+\Delta t}$ is not known at time $t$. Indeed, in the application considered here, the flow of water $Q$ and electricity price $P$ are stochastic. Therefore, we do not have sufficient information to determine the optimal choice in \eqref{eq:optimal} or the value in \eqref{eq:optimalvalue}.

To remedy this lack of information we create \textit{approximately} optimal strategies based on historical and forecast estimates using \eqref{eq:strategy}. 
To be more precise, let $\bar  X=(\bar Q, \bar H, \bar P)$ denote a known deterministic estimate of the underlying processes, possibly respecting forecasts. At each time $t$, we now proceed as outlined above, with the difference that we replace the unknown stochastic variable $X^x_{t+\Delta t}$ with its \textit{deterministic} counterpart $\bar X^x_{{t+\Delta t}}$. Given a terminal value $\{\bar V_1 (T,x),\dots,\bar V_m(T,x)\}$ we can thus recursively construct an \textit{approximate} value function $\{\bar V_1(t,x), \dots, \bar V_m(t,x)\}$ by mimicking \eqref{eq:optimalvalue}, i.e., 
\begin{eqnarray*}
&&\bar V_i(t, x) = \max_{j \in \{1,2,\dots,m\}}  \left\{\phi_j(t,x) \cdot \Delta t+\bar V_j(t+\Delta t, \bar X^x_{t+\Delta t}) - c_{ij}(t, x)\right\}, \notag \\ \qquad && \bar V_i(T,x) = g_i(x).
\end{eqnarray*}
The function $\bar V(t,x)$ does not coincide with $V(t,x)$ in our original problem as it is based merely on estimates on $X_t$. However, one suspects that moving from mode $i$ to mode $j$ whenever
$$
\bar V_i(t,X_t) = \bar V_j (t,X_t) - c_{ij}(t,X_t)
$$
should be close to optimal in the original optimization problem, regardless of the actual value of the function $\bar V$. Indeed, we expect the stochastic process $X_t$ to behave in a similar fashion to $\bar X$, so the optimal strategy should not be too different either. In particular, if the process $\bar X$ includes short-term forecasts, the decision in the short run should be close to optimal since the estimate of the nearby future is then very good and the long term effects should be respected by the estimate $\bar X$. 

\subsection{River flow model}
\label{sec:river-model}

We will denote our modeled flow in the river at time $s$ by $\bar Q_s$ and we let $\bar q = \bar q(s)$ be an historical estimate of the flow at that time. 
In particular, we let $\bar q$ be a 7-day standard moving average based on data from 1980-2014. 
In case of no (perfect) forecast, $\bar Q_s$ is known only up to the present time $t$ whereas the historical estimate $\bar q$ is known for the entire year.

In case no forecast is available, we assume that the modeled flow $\bar Q_s$, $s>t$, reverts towards $\bar q$ so that the difference, $\bar\lambda(s) = \bar Q_s - \bar q(s)$, satisfies

$$
\frac{d \bar\lambda}{ds} = -\kappa \bar\lambda(s) 
\quad \text{i.e. }\quad \bar\lambda(s) = \bar\lambda({t}) e^{-\kappa (s - t)}
\quad \text{for} \quad s \geq t.
$$
Thus 
\begin{align*}
\bar Q_s = \left(\bar Q_t - \bar q(t)\right) e^{-\kappa (s - t)} +\bar q(s) 
\quad \text{for} \quad s \geq t.
\end{align*}
This means that an initial difference from the historical flow $q$ vanishes exponentially. 
Writing this in terms of the half life $T_{1/2}$ we have
\begin{align}\label{eq:flowmodel_new}
\bar Q_s = \left(\bar Q_t -\bar q(t)\right) 2^{-\frac{s-t}{T_{1/2}}} +\bar q(s) 
\quad \text{for} \quad s \geq t.
\end{align}

The flow data used in our numerical examination is from Sävarån in the northern parts of Sweden and is gathered from the Swedish Meteorological and Hydrological Institute.\footnote{Flow data was downloaded from http://vattenwebb.smhi.se/station (station number 2236) on September 1, 2023} 
The average flow $\bar q$ is based on data from 1980-2014 while the data from 2015-2022 is used solely for testing our optimization method in Section \ref{sec:results}. 
The mean flow $\bar q$ together with the flow of the bench-marking years are shown in Figure \ref{fig:flow}. 
We set the half time without further investigations to $T_{1/2}=10$ days and investigate the sensitivity of our results due to this choice later on, see Remark \ref{remark:OSPflow} in Section \ref{sec:results}.
A direction field of our simple flow model in \eqref{eq:flowmodel_new} is shown in Figure \ref{fig:flowfield}.

\begin{figure}
    \centering
    \includegraphics[scale=0.2]{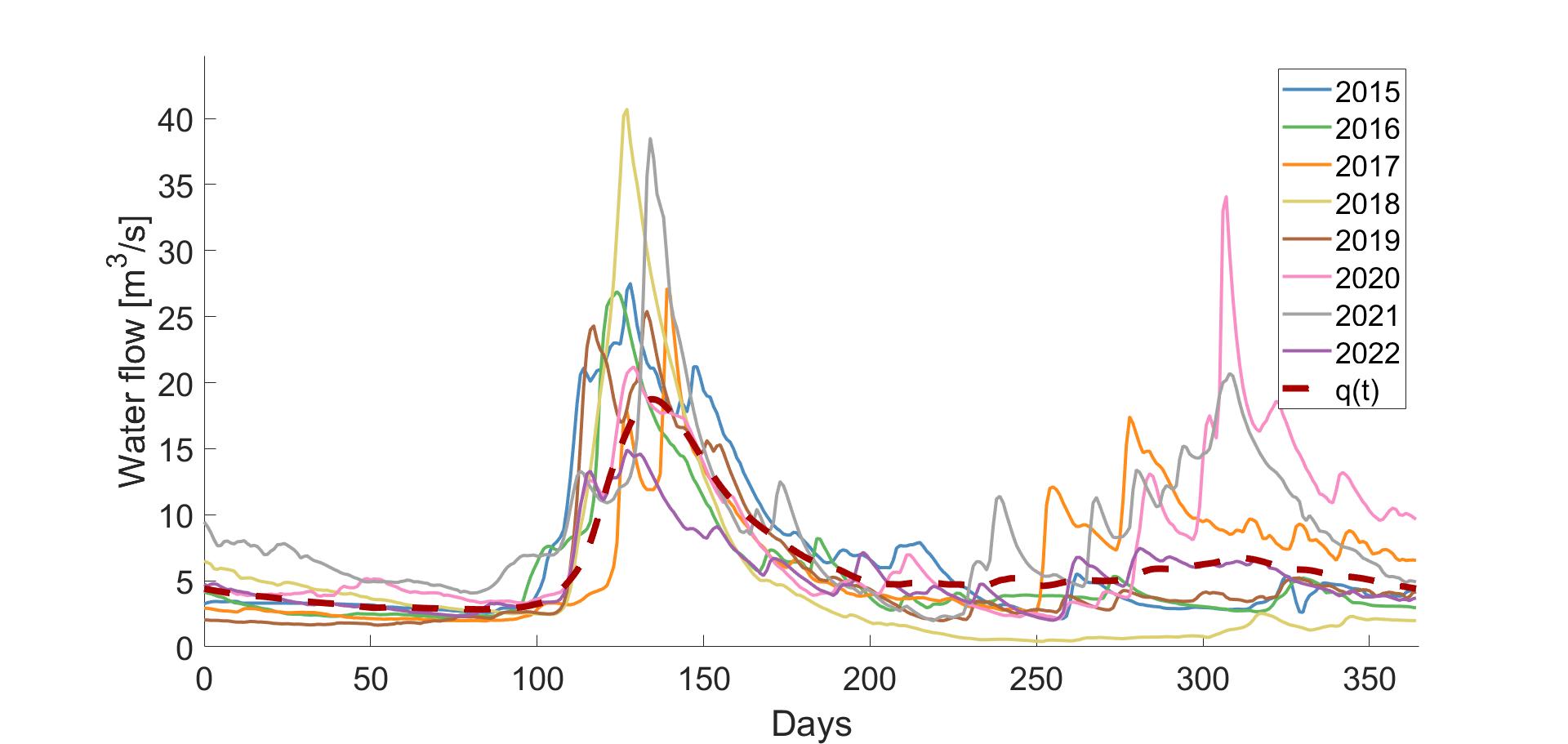}
    \caption{Mean flow $\bar q$ (dashed) based on data from 1980-2014 together with actual flows from 2015-2022. Flow data is from Sävarån in the northern parts of Sweden.}
    \label{fig:flow}
\end{figure}

\begin{figure}
    \includegraphics[scale=0.25]{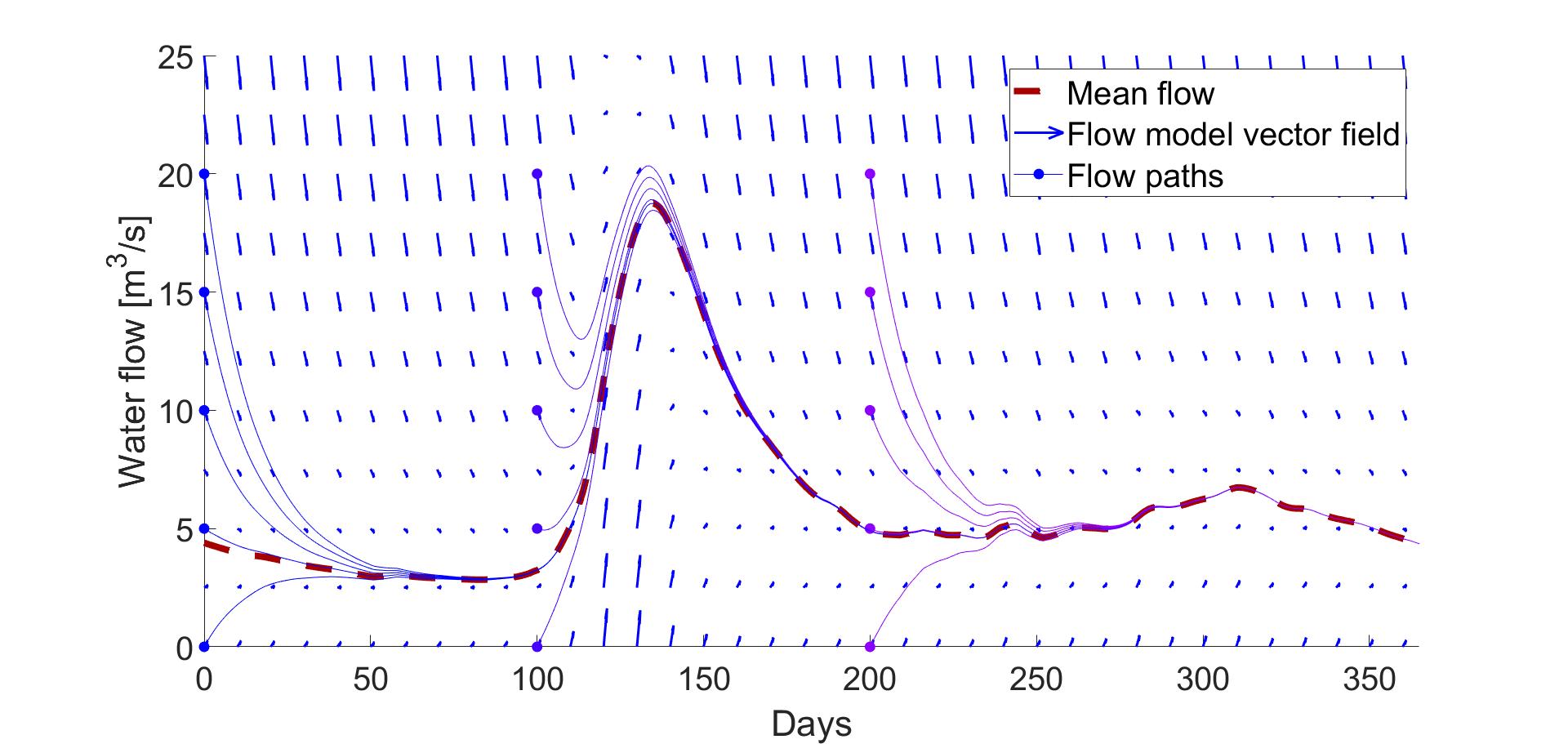}
    \caption{The flow field of our model, showing how the model flow reverts towards the mean flow $\bar q$ with $T_{1/2}=10$ days.}
    \label{fig:flowfield}
\end{figure}

When a forecast of $M$ days is available, we replace the first $M$ days of $\bar Q_s$ with the corresponding forecast. After these days, the modeled flow $\bar Q_s$ is assumed to return to the mean flow $\bar q$ exponentially as outlined above. To avoid forecast modeling, which is not the topic of the current paper, we simply assume the forecasts are perfect and use the actual flow as prediction in our numerical investigation below. We briefly comment on the impact of the specific flow model constructed here in Remark \ref{remark:OSPflow}.

\subsection{Power plant modeling}

When considering hydropower plants with a dam, it is natural to model the power output of each unit, i.e., turbine and generator pair, as a function of the head $H_t$ of the reservoir and the flow of water $F_t$ through the turbine. 
We thus assume that the payoff from the power plant depends on the controlled processes $H_t = H^{\alpha}_t$ and $F_t = F^{\alpha}_t$ 
in which $\alpha$ is the control. The head is given by
\begin{equation*} 
\frac{dH^\alpha}{dt}= g_H(F_t^\alpha, Q_t,H^\alpha_t), \qquad H_0=h
\end{equation*}
where $Q_t$ is the inflow to the reservoir (i.e. the river flow as above), $\alpha$ indicates the current production mode, and $g_H$ is a function given by the shape of the reservoir. Note in particular that the chosen strategy $\alpha$ has a direct impact on the dynamics of the water head $H^\alpha_t$. 
For the sake of our numerical example, we assume that the dam has the simple shape of a cone with maximum height $H_{max}$ and that it can hold enough water to supply the power plant with water for its design speed $F_d$ $m^3/s$ during $N$ days (see \eqref{eq:efficiency} for a definition of design speed). 
From here on in, we will refer to the size of the dam in terms of this number $N$ and in our investigation we will look at a wide variety of this value.
{Roughly speaking, a large (respectively small) value for N means that we consider a power plant with a dam that is large (small) relative to the size of its turbine.}
Simple arithmetic gives that
$$
H^\alpha_t= H_{max} \left (\frac{V^\alpha_t}{V_{max}} \right ) ^{1/3}
$$  
where $V^\alpha_t$ is the amount of water in the reservoir at time $t$, and $H_{max}$ and $V_{max}$ are the height and capacity of the reservoir, respectively. 

We assume that the plant consists of a single unit that can generate electricity for all flows between $F_{min}$ and $F_{max}$. We normalize all data so that the payoff when production is completely shut down ($i=0$) is $0$. When in productive mode, we follow the reasoning in \cite[section 4]{OOL22} and let the payoff be given by
\begin{align}\label{eq:f2dam}
\phi(F_t^\alpha,H_t,P_t) = - c_{run} + \begin{cases}
-c_{low} & \mbox{if $H_t =0$} \\
    \rho g H_t\, \eta(F^\alpha_t)\,F^\alpha_t \, P_t 		& \mbox{if $H_t >0$},	
\end{cases} 
\end{align}
where $F_t^\alpha$ is the amount of water run through the generator, $\rho=10^3 kg/m^3$, $g=9.82 m/s^2$, $c_{run}$ and $c_{low}$ are constants, and $\eta$ is an efficiency curve 
\begin{equation} \label{eq:efficiency}
\eta(F) = \hat \alpha - \hat \beta \left( \frac{F}{F_d} - 1\right)^2
\end{equation}
specific for the unit under consideration. 
{Here, $F_d$ is the design speed and $\hat \alpha$ and $\hat \beta$ are constants, see Table \ref{tab:parameters}, Figure \ref{fig:efficiency} and \cite[Figure 2]{OOL22}; the efficiency curve comes from fitting data for a Swedish Kaplan turbine.} The condition $H_t =0$ in \eqref{eq:f2dam} implies a penalty if the dam runs empty, thereby failing to meet the minimum requirements of the unit. Note that the running cost may exceed the possible profit from generating electricity if the water head is too low, so the dam may be effectively ''empty" before $H_t=0$ (but production is nevertheless possible without penalization as long as $H_t>0$).

\begin{figure}
    \centering
    \includegraphics[scale=0.2]{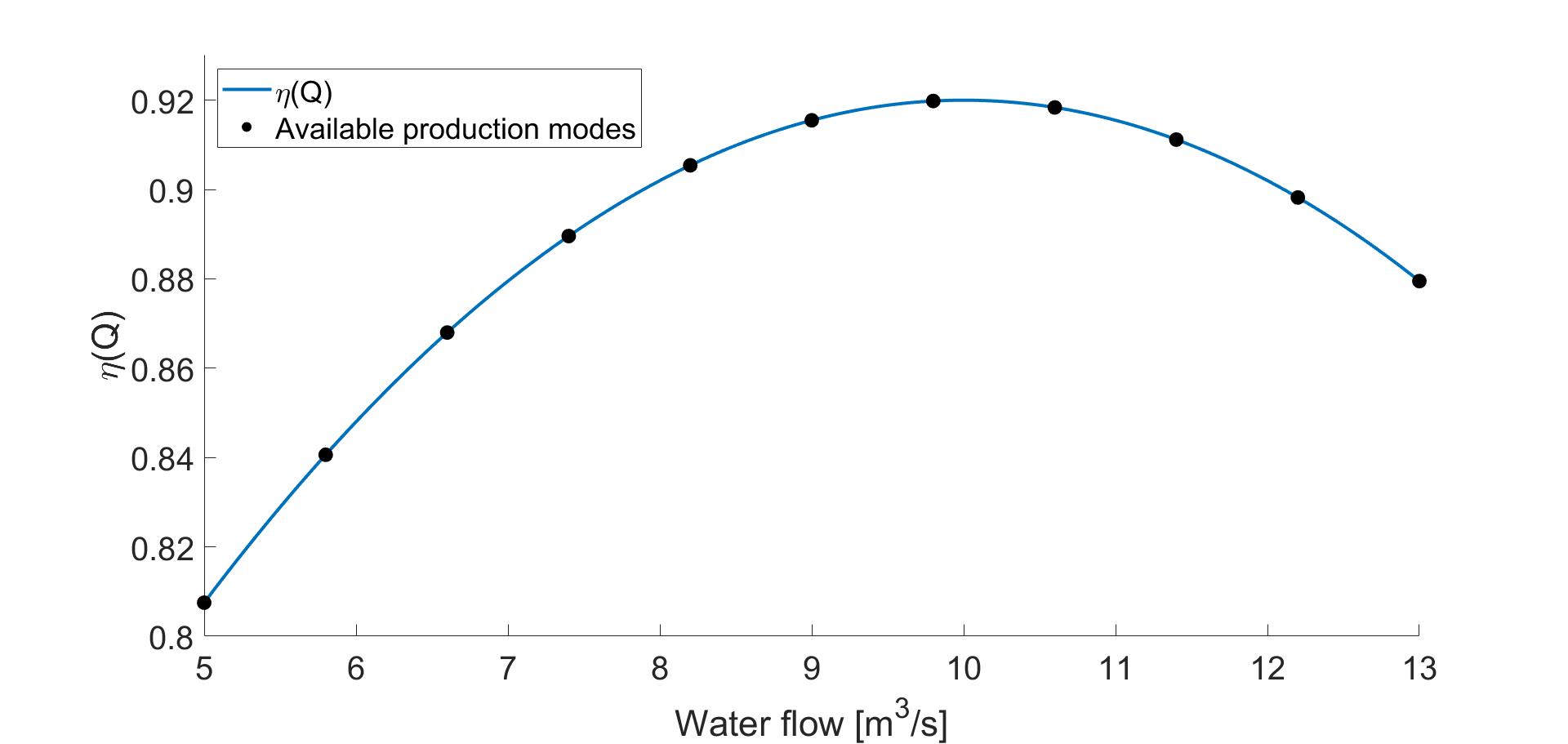}
    \caption{The efficiency curve given by \eqref{eq:efficiency}, with $\hat \alpha = 0.92, \hat \beta = 0.45$, together with the available production modes of the unit.}
    \label{fig:efficiency}
\end{figure}
   
\begin{table} 
\centering
   \begin{tabular}{| l | l | l | l | }
    \hline
$P_0$        & 1 $[m.u. /kWh]$ &$c_{low}$      & $ 1000 $ $[m.u./h]$ \\  \hline
$H_{max}$    & 5  $[m]$        & $c_{run}$      & $ 100 $ $[m.u./h]$ \\ \hline
T           & 365 [days]      & $F_{min}$      & 5 $[m^3/s]$  \\ \hline
$\hat\alpha$    & 0.92            &$F_{max}$      &13 $[m^3/s]$   \\ \hline
$\hat\beta$     & 0.45            & $F_{d}$        &10 $[m^3/s]$    \\ \hline
$T_{1/2}$   & 10 [days]           &      &    \\ \hline
\end{tabular}
\caption{Parameter values used in our numerical investigation.}
\label{tab:parameters}
\end{table}%

{The turbine can be adjusted to fit the current flow of water but each such adjustment is associated with a cost.
This cost, called a ``switching cost", is set as a fraction of the profit generated by the plant if it works at maximum capacity for a full year at unit electricity price and without interruptions. In this particular case, this maximum is given by 
\begin{equation}\label{eq:D}
D=\phi(F_{max},H_{max},1)\cdot 365=\phi_{9}(H_{max},1)\cdot 365
\end{equation}
and the cost of adjustment is $\gamma D$, where $\gamma \in [0, 0.01]$ is a parameter showing, in general, how costly it is to make this type of change in the production. 
Furthermore, we assume the cost of starting/stopping the generator is $25$ times that of simply adjusting an already running generator, see \eqref{eq:sw-costs-def} below. 

\subsection{Further assumptions and numerical procedure}

Our optimization is based on the time discretization $[0:\Delta t:T]$ and spatial discretization
\begin{equation}\label{eq:discretizations}
Q=[0:\Delta q: \hat Q], \qquad H=[0:\Delta h: H_{max} ]
\end{equation}
where $\Delta t= 1$ day, $\Delta q =1/4$ $m^3/s$, $\hat Q=2 \cdot F_{max}$, and $\Delta h$ corresponds to $\frac{1}{4000}$ of the total dam size. For computational ease, all quantities are calculated on a grid point of the discretization \eqref{eq:discretizations} using rounding to the nearest point when necessary. A finer (or coarser) grid can therefore alter the payoffs and corresponding strategies slightly but not enough to change the qualitative results.

To capture the natural seasonality of the problem we consider an optimization horizon of $T=365$ days\footnote{Leap days are excluded for simplicity of presentation.} and allow the manager to change the mode of production once per day. The electricity price is taken to be constant $P_t \equiv P_0 =1$, corresponding roughly to maximizing the output of electricity rather than the monetary profit\footnote{Time-dependent electricity price can be handled without any additional complications but requires slightly longer computational time and obstructs the interpretation of the results.}.  
For simplicity of presentation, we require the plant to end up in the same production mode as it started in (i.e., $i=0$, ''off"). 

 For our numerical calculations, we restrict the adjustment of the turbine to 9 different modes ($i\in \{1,2,\dots, 9\}$), giving in total 10 different modes of production, $i=0$ meaning no production and the remaining modes having $F^\alpha$ spread evenly from $F_{min}$ to $F_{max}$,
$$
F_i = F_{min} + \frac{(i-1)}{8}(F_{max}-F_{min}) \quad \text{for} \quad i \in \{1,2, \dots,9\}.
$$
For the corresponding running payoff, cf. \eqref{eq:f2dam}, we use the notation 
$$
\phi_0=0 \qquad \mbox{and} \qquad \phi_i (H_t,P_t) := \phi(F_i,H_t,P_t).
$$
(This notation is already used in \eqref{eq:D} above.) Similarly, the cost of switching between the different production modes $i \in \{0, \dots, 9\}$ is given by
\begin{align}\label{eq:sw-costs-def}
c_{ij} =  \begin{cases}
    0 & \mbox{if $i=j$} \\
    \gamma D  & \mbox{if $i \neq j$ and ($i=0$ or $j=0$)} \\
    \frac{\gamma D}{25}  & \mbox{if $i \neq j$ and $i,j \neq0$} ,
\end{cases} 
\end{align}
The efficiency curve $\eta(F)$ used is depicted in Figure \ref{fig:efficiency} with the corresponding allowed production flows $F_i$ marked with dots. 
}
As water in the reservoir has value we must take any change in the reservoir from the beginning to the end of the optimization period into account in the final result. This is done by establishing the value of water in the reservoir as the profit this water would generate if used to run the generator at design speed $F_d$, disregarding running costs. 
Any change in the reservoir from the initial level $H_0 =H_{max}$ (which corresponds to the dam being full) is adjusted in the final profit. Note that the assumption of design speed and no running cost implies a larger penalty for missing water than what could be gained from using it for production, thereby forcing an optimal strategy to end with the dam full.

Naturally, the numerical values to be used vary with the specific problem, river, and power plant under consideration. The parameter values applied here are summarized and presented in Table \ref{tab:parameters}; we refer to \cite{OOL22} for details and motivations. The exact values should have little impact on the qualitative nature of our results. For our model and the sake of this paper, the most significant parameters are the forecast length, dam size, and switching cost. If nothing else is specified we consider $M=10$ days forecast, a dam size of $N=30$ days at design speed, and switching cost parameter $\gamma=0.0025$. In the next section, we vary these parameters one-by-one to highlight their impact on the optimal strategies and the end result. 

\section{Results}
\label{sec:results}

The suggested production schemes perform remarkably well, averaging $97.5$~\% of the theoretical maximum over the years 2015-2022 for parameters as above (see Figures \ref{fig:forecastdamsize} and \ref{fig:resultsswcost} for the performance with other parameters).
A detailed view of the optimal and suggested production scheme for 2022 is presented in Figure \ref{fig:strategydetailed}, in which we 
observe the following: 
The strategies typically avoid running the plant at full capacity (which is natural because of the efficiency of the unit, recall Figure \ref{fig:efficiency}) 
and keep the the reservoir head above $80 \%$ (which corresponds to about $60 \%$ of the dam capacity) and at about $90\%$ on average over the year. 
Concerning differences and similarities between the DPP strategy and the optimal strategy, we observe that the main differences occur during and after the spring flood, which is likely due to the fact that the flow fluctuates more during this period. 
\begin{landscape}
\begin{figure}
\centering
   \includegraphics[width=1.5\textwidth]{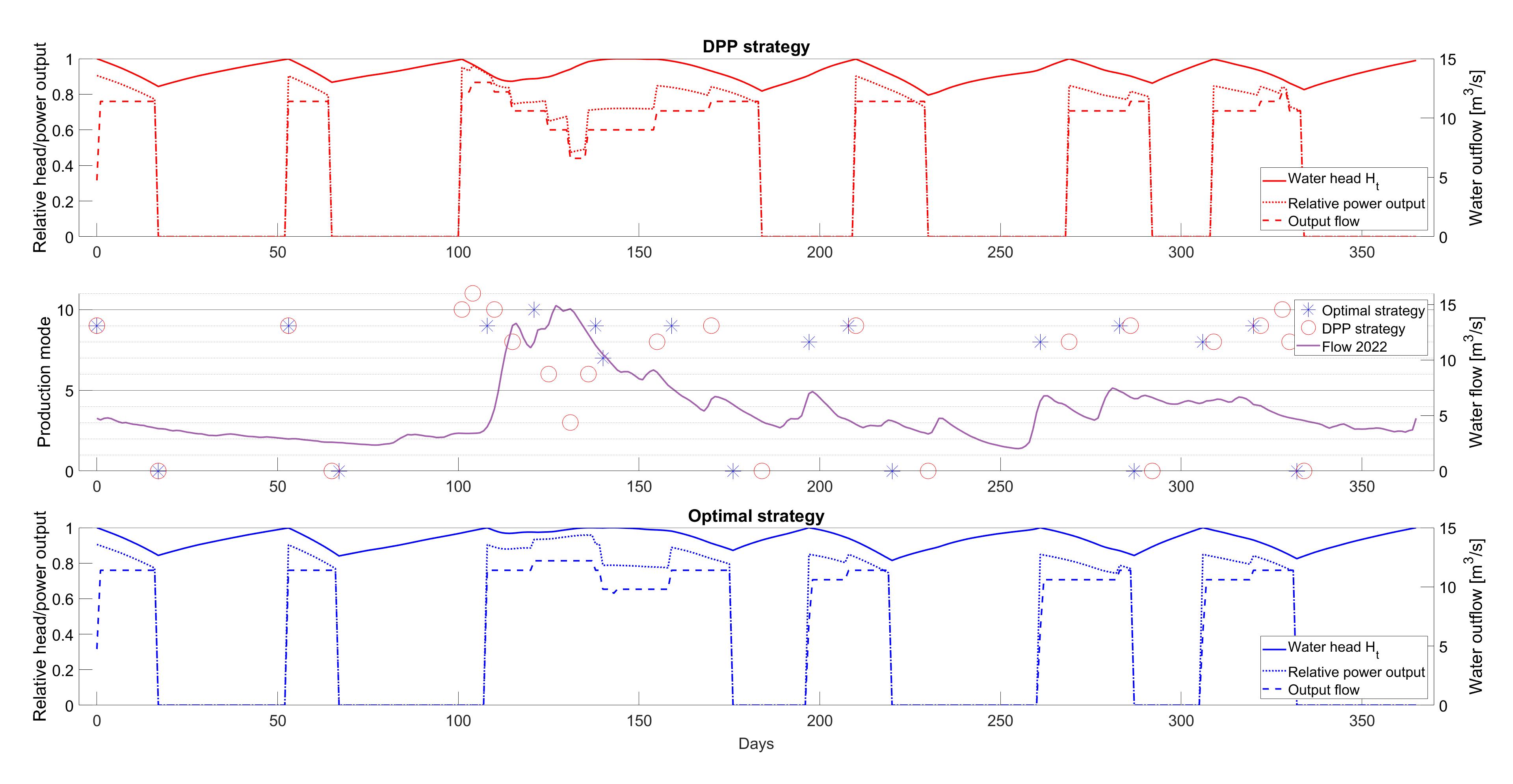}   
\caption{Production strategies for year $2022$ with a dam size of $N=30$ days, $\gamma=0.0025$ and forecast length $10$ days. The suggested strategy (red) gives $96.5$ \% of that of the optimal strategy (blue).
The presence of a marker (blue star for the optimal solution and red circle for the DPP solution) indicates a change in production, with the height of the marker (left axis) representing the mode held until the next change.}
\label{fig:strategydetailed}
\end{figure}
\end{landscape}

\subsection{Dam size}

The presence of a dam significantly increases the management options and therefore also the payoff of the power plant.
Another benefit of a dam is that it reduces the importance of accurate short-term forecasts, as can be seen from Figure \ref{fig:forecastdamsize}. This is in line with what should be expected as the storage of water can be used to manually counteract sudden changes in the river flow to keep production efficient. Our strategy has no difficulties finding these adjustments. The need for a short-term forecast vanishes as the dam grows as the current inflow then becomes insignificant in comparison with the long-term average.
With a small dam one must use a wider range of operating modes and turn the plant on/off more often whereas with a large dam, a fairly decent result can be achieved using a single mode of operation, see Figure~\ref{fig:optstratvaryingdam2022}.
Our method performs well in both cases.

\begin{figure}
    \centering
        \includegraphics[width=\textwidth]{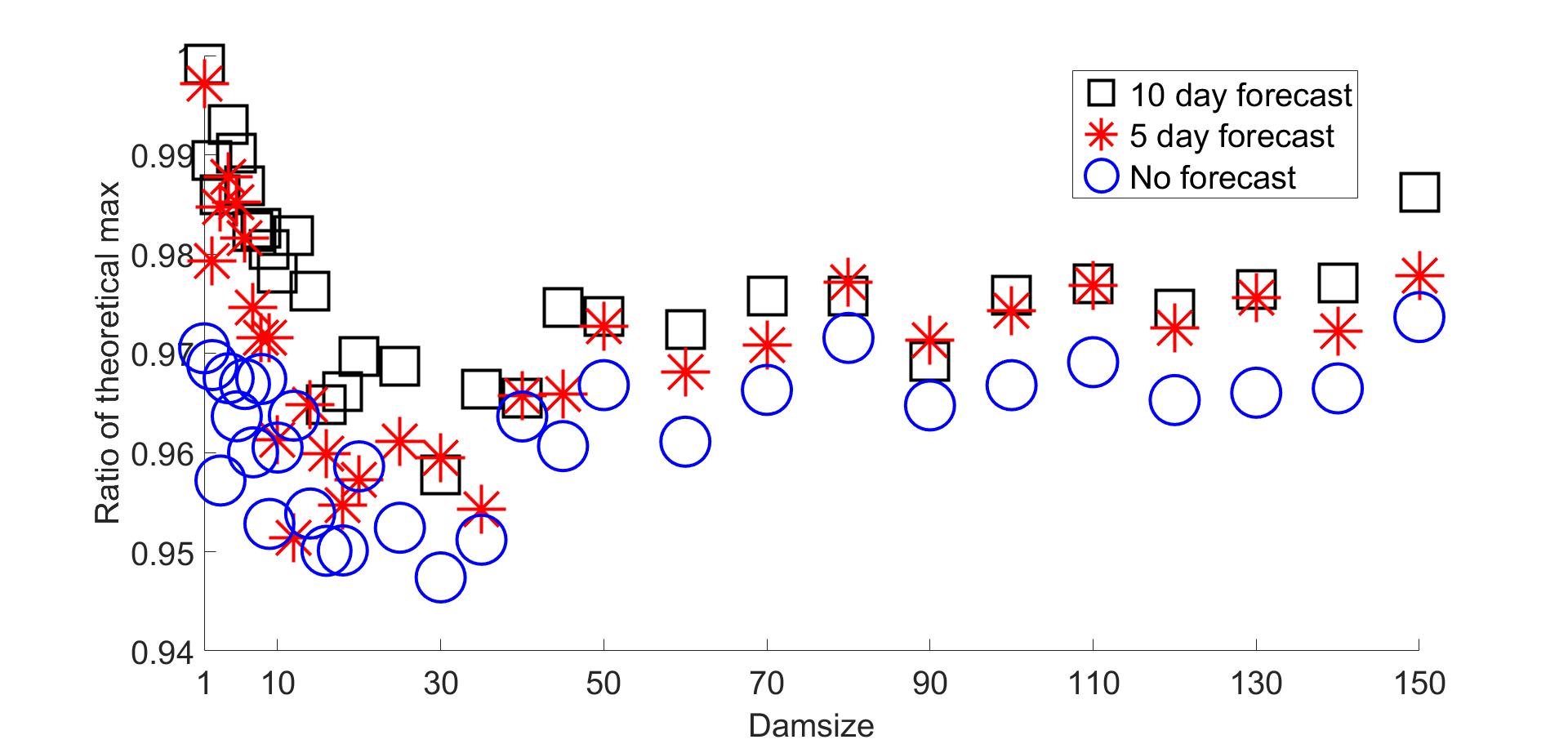}
    \caption{The benefits of a forecast become less explicit as the size of the dam grows. Results are the average over the years 2015-2022 with $\gamma=0.0025$.
    }
    \label{fig:forecastdamsize}
\end{figure}

\begin{figure}
\centering
    \includegraphics[width=\textwidth]{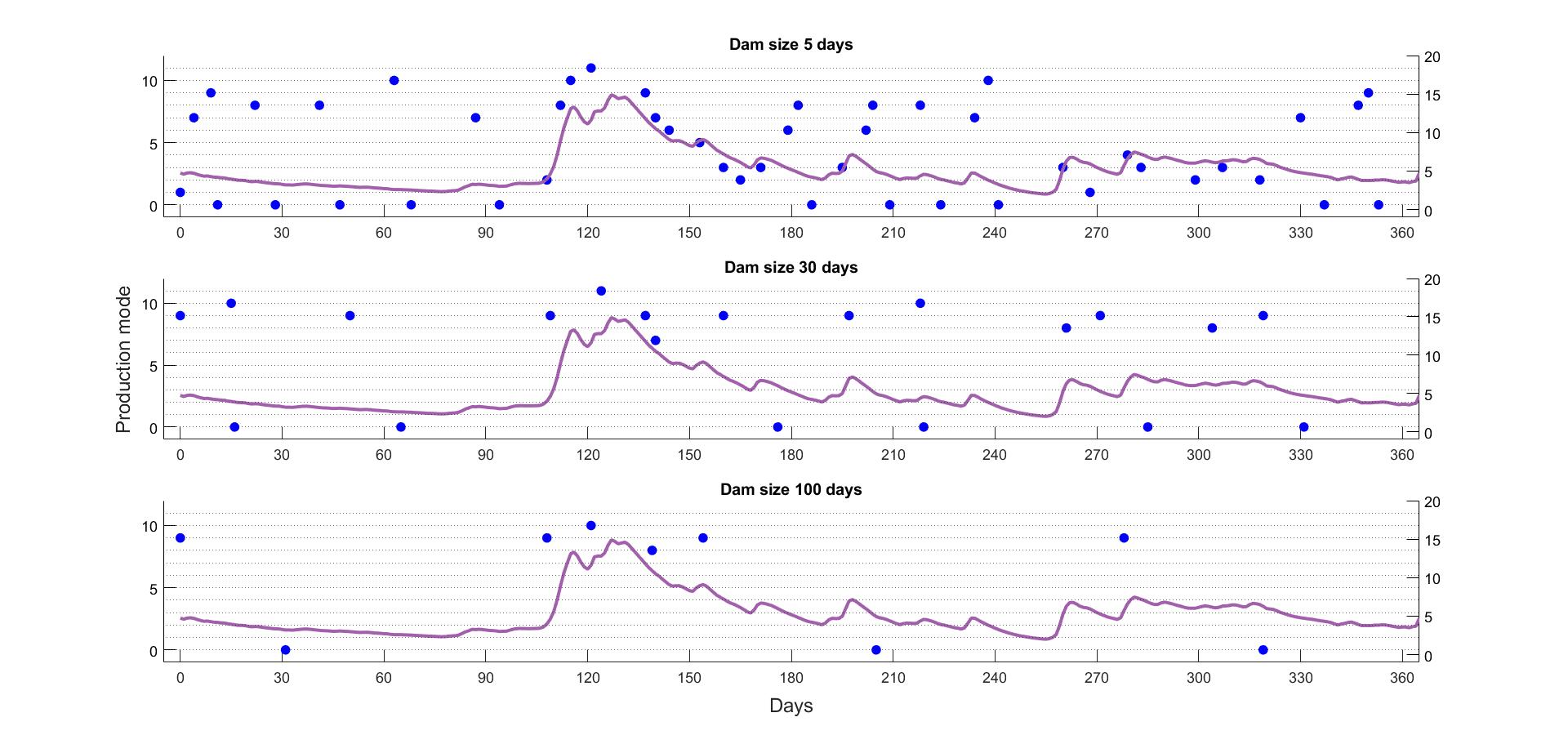}
    \caption{Optimal strategy for dam size $N=\{5,30,100\}$ days. The need for active management increases with a dam but vanishes as the dam grows. The presence of a marker (blue dot) indicates a change in production, with the height of the marker (left axis) representing the mode held until next change. Flow data for year 2022 with $\gamma=0.0025$. }
    \label{fig:optstratvaryingdam2022}
\end{figure}

\subsection{Switching cost}
Clearly, the total profit decreases as the cost of changing production mode increases. As with the dam size, the cost of changing production also affects how many switches should be made in an optimal strategy; smaller costs lead to more active management and vice versa, see Figure \ref{fig:optstratvaryingswcost2022}. We also observe that the performance decreases with increasing (but still reasonable) switching costs, which likely is due to the higher cost of reversing a non-optimal decision.
Our scheme is stable w.r.t. these changes in the sense that it adjusts the suggested strategy as necessary to find an efficient production plan in all cases, see Figure \ref{fig:resultsswcost}. 

\begin{figure}
\centering
    \includegraphics[width=\textwidth]{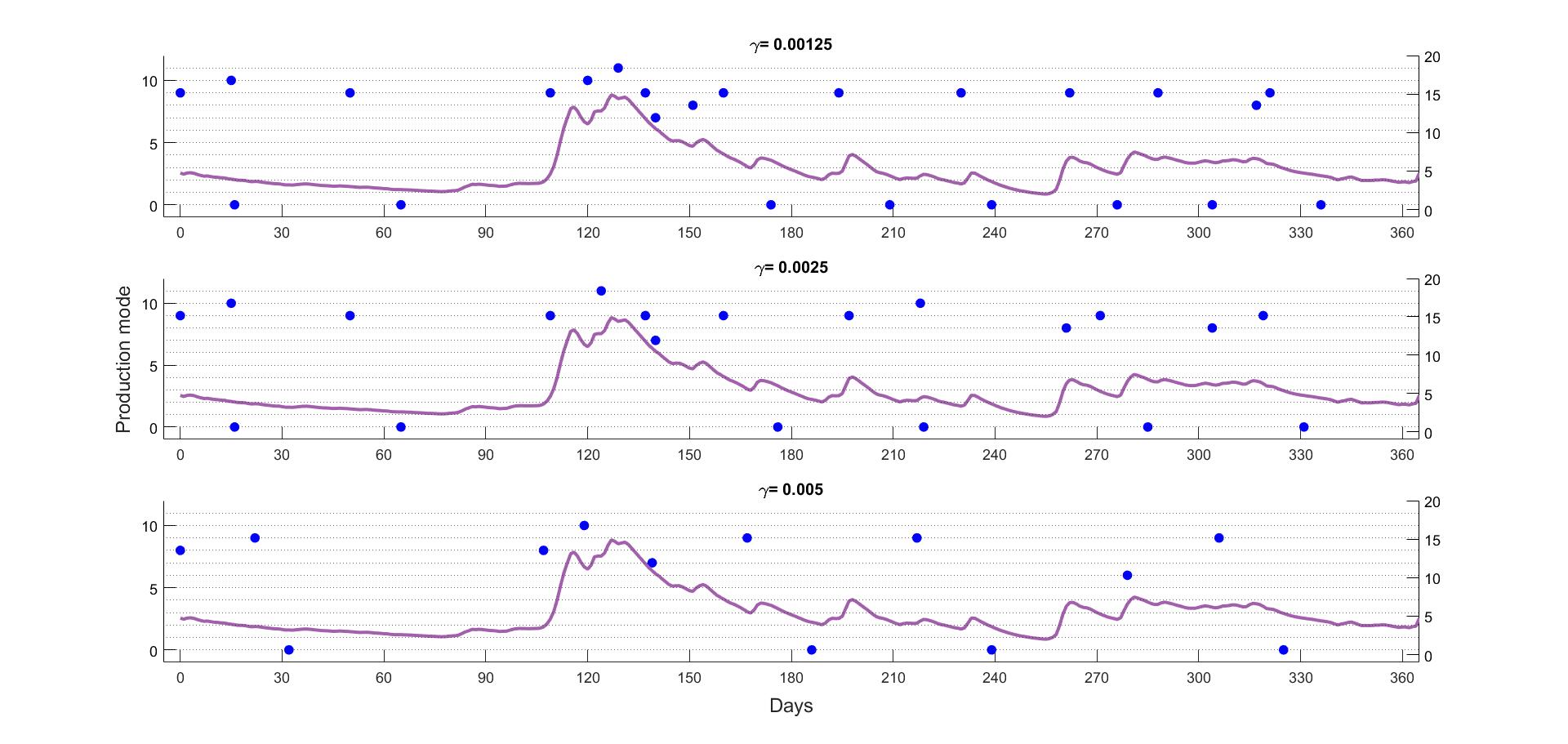}
    \caption{Optimal strategies for $\gamma=\{0.00125, 0.0025, 0.005\}$. Lower switching cost naturally leads to more changes in the production. The presence of a marker (blue dot) indicates a change in production, with the height of the marker (left axis) representing the mode held until next change. Flow from year 2022 with dam size $N=30$ days.}
    \label{fig:optstratvaryingswcost2022}
\end{figure}

\begin{figure}
    \centering
    \includegraphics[scale=.2]{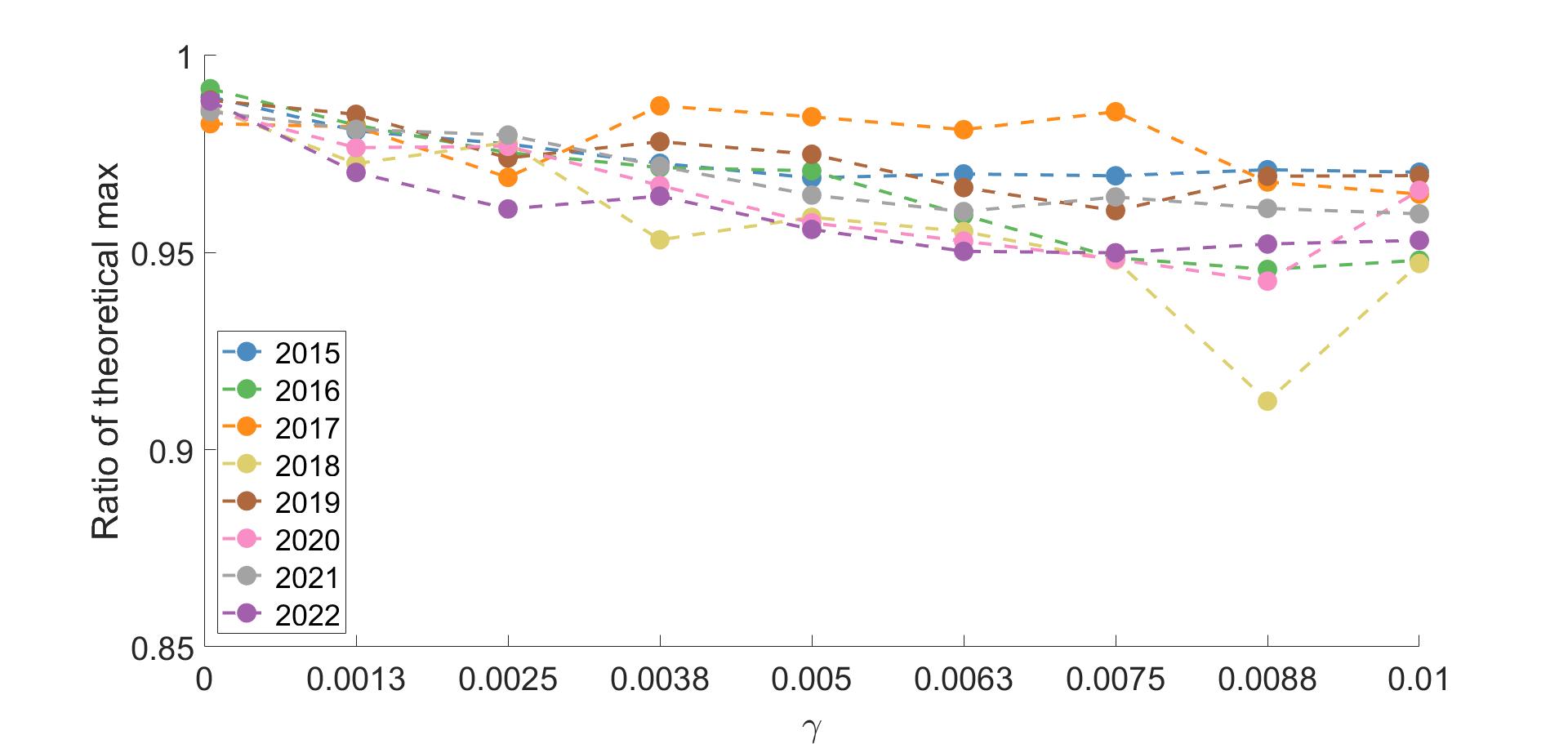}
    \caption{We obtain results close to the theoretical maximum for a wide range of switching costs, i.e., different values of $\gamma$. Our method thus appears stable w.r.t. to this parameter.}
    \label{fig:resultsswcost}
\end{figure}

\section{Managing a run of river power plant - comparison to optimal switching}
\label{sec:oldsetup}

The objective function \eqref{eq:tooptimize} investigated here falls into the framework of optimal switching theory. This theory was used in \cite{OOL22} for production planning of a Run-of-River (RoR) power plant with two units that could be regulated and switched on and off depending on the natural flow of the river. In essence, this corresponds to having three different modes of production and no dam to store water. For comparison of methods, we here mimic that setup and apply our much simpler method to the same data set. We present and compare the results of the method outlined in Section \ref{sec:method} and that of \cite{OOL22} (named OSP below) for the years 2019-2022. Throughout this section, we refer to these methods as DPP and OSP, respectively. 
The parameters are as in Section \ref{sec:results} and \cite{OOL22}, as applicable.

The plant can be run in three different modes; shut down (mode $i=0$), 1 unit running ($i=1$), or 2 units running ($i=2$). Both generators have the efficiency $\eta$ as in Figure \ref{fig:efficiency}. When in mode $i=2$, the now uncontrolled flow of water can be split between the two units at no cost to maximize the combined output of the pair. The data is normalized so that the payoff from mode $0$ (shut down) is $0$, i.e., $\phi_0\equiv 0$ and in productive mode the payoff is given by
\begin{align*}
\phi_1(F_t,P_t) &= - c_{run} +
\begin{cases}
-c_{low}				& \textrm{if}\quad F_t < F_{min},\\
c\, \eta(F_t)\,F_t \,P_t			& \textrm{if}\quad F_{min} \leq F_t < F_{max},\\
c\, \eta(F_{max}) F_{max} \,P_t			& \textrm{if}\quad F_{max} \leq F_t,
\end{cases} \\
\phi_2(F_t,P_t)&= \max_{\delta \in [0,1]} \left\{ \phi_1 (\delta F_t,P_t) + \phi_1 ((1-\delta) F_t,P_t) \right \},
\end{align*}
where $c=\rho g H_t$. Note that the water head $H_t\equiv H_{max}$ is fixed and $F_t=Q_t$ as the plant cannot store any water. As above, $P_t\equiv 1$ and the cost of switching is defined via \eqref{eq:D}
as
\begin{equation*}
c_{ij}=\begin{cases}
    0       & \mbox{if $|i=j|=0$}    \\
    \gamma D       & \mbox{if $|i-j|=1$}\\
   1.5\cdot \gamma D       & \mbox{if $|i-j|=2$}
\end{cases}.
\end{equation*}

As indicated already by the results in Figure \ref{fig:forecastdamsize}, accurate forecasts are of great importance for RoR power plants when sudden changes of the water flow cannot be counteracted by stored water. However, the efficiency of the OSP method is not as sensitive to forecasts as the DPP method above, see Figures \ref{fig:difference1} and \ref{fig:difference2}. This is due to the inherent stochastic features of the OSP method which leads to more "wait-and-see" strategies than the DPP method.  
The OSP method is stable but rarely finds the true optimal strategy, i.e., the performance ratio (ratio of the payoff to the theoretical maximal payoff) is $<1$, see Figure \ref{fig:difference1}, while the DPP method in many cases finds the true optimum. On the other hand, OSP avoids major pitfalls and performs close to optimal on most occasions, both with and without forecast, while the DPP method is more prone to make costly sub-optimal decisions, especially when the information on future flow is limited, see Figure \ref{fig:difference2}. We consider such a specific scenario below to give an intuition of how the schemes differ in their decision-making.

\begin{figure}
\centering
   \begin{subfigure}[a]{1\textwidth}
   \includegraphics[width=\textwidth]{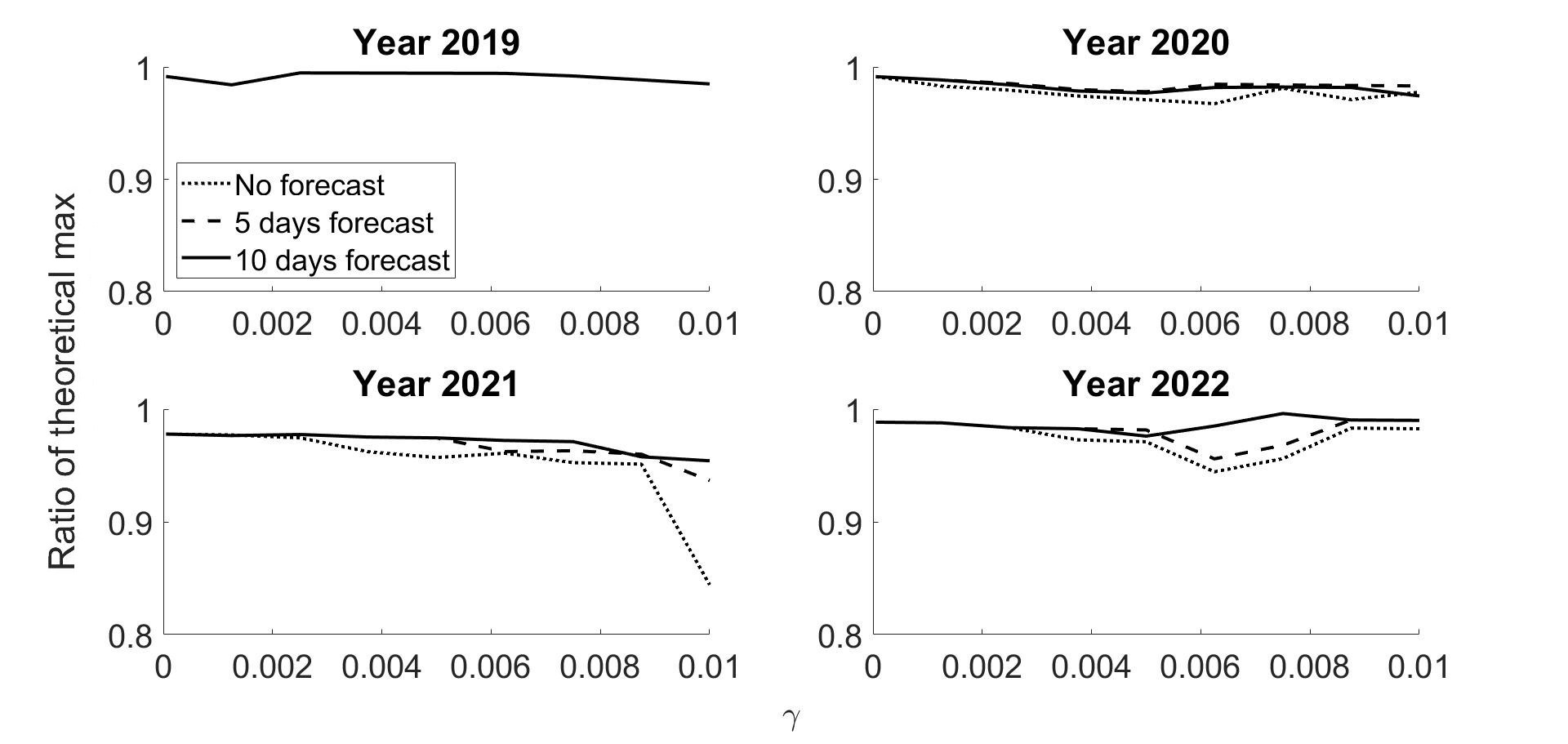}
    \caption{Results for the OSP-method of \cite{OOL22}.}
    \label{fig:difference1}
\end{subfigure}

\begin{subfigure}[b]{1\textwidth}
    \includegraphics[width=\textwidth]{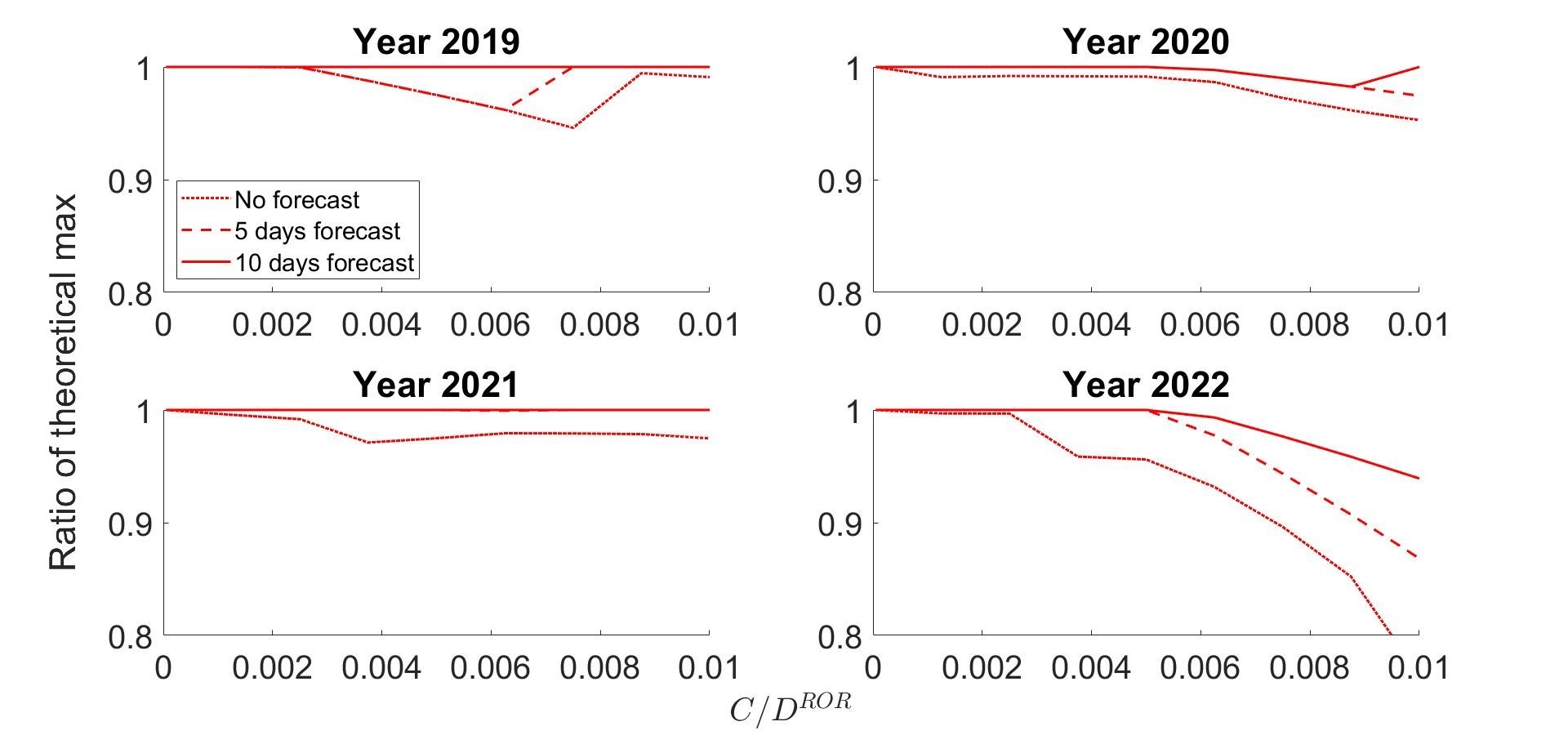}
    \caption{The DPP-method using the simple flow model \eqref{eq:flowmodel_new}.}
    \label{fig:difference2}
\end{subfigure}
\caption{The strategies developed in this paper often perform better than the more complicated approach of \cite{OOL22} when the forecast is sufficiently good. However, when information is scarce or incorrect, the stochastic model of \cite{OOL22} typically outperforms our method. Note in particular the difference for year 2022 where spring flood differs much from the historical average, cf. Figure~\ref{fig:flow}.}
\end{figure}

For year 2019 and parameters $\gamma=0.0075$, and $M=5$ days the performance ratio (ratio of the payoff to the theoretical maximal payoff) is $0.992$ for the OSP strategy and $1.000$ for the DPP strategy. That is, the DPP strategy finds the true optimum while the OSP strategy comes close, but not all the way. The corresponding strategies are shown in Figure \ref{fig:stratdiff2019} and we scrutinize two of these decisions to show where the performance difference appears.
At time $t=110$, the OSP strategy chooses to open production at the intermediate level ($i=1$) although the optimal strategy is to wait and open at full capacity later on at $t=112$. The situation at the time of decision is depicted in Figure \ref{fig:closerinspection2019a}. The OSP model forecasts a flow below that of DPP and, in addition, anticipates deviations from this forecast, resulting in the safer choice of an intermediate step in the production. Note however that the OSP takes this action \textit{before} the optimal strategy opens to mode $i=2$, so some of the loss due to the extra switch is regained. 

Figure \ref{fig:closerinspection2019b} examines the point $t=325$ where both strategies optimally refrain from turning production on despite going into the profitable region, $F=Q>F_{min}$. The forecast is sufficiently long for the DPP model to detect the upcoming decrease in flow and avoid turning production on. The reason for the OSP model to refrain from switching mode is different; it projects a flow \textit{above} the critical level $Q_{min}$ for a sufficiently long time for a change in production(back and forth) to be profitable but anticipates deviations from this forecast and therefore requires more margins before taking action. With a slightly shorter forecast of 4 days, the projected flow is sufficiently above the critical level for the DPP strategy to sub-optimally turn production on at $t=325$ (and then off again at $t=333$) while the OSP strategy with its built-in stochastic features remains unchanged and optimal even with this shorter forecast.
     
    \begin{figure}
    \centering
    \includegraphics[width=\textwidth]{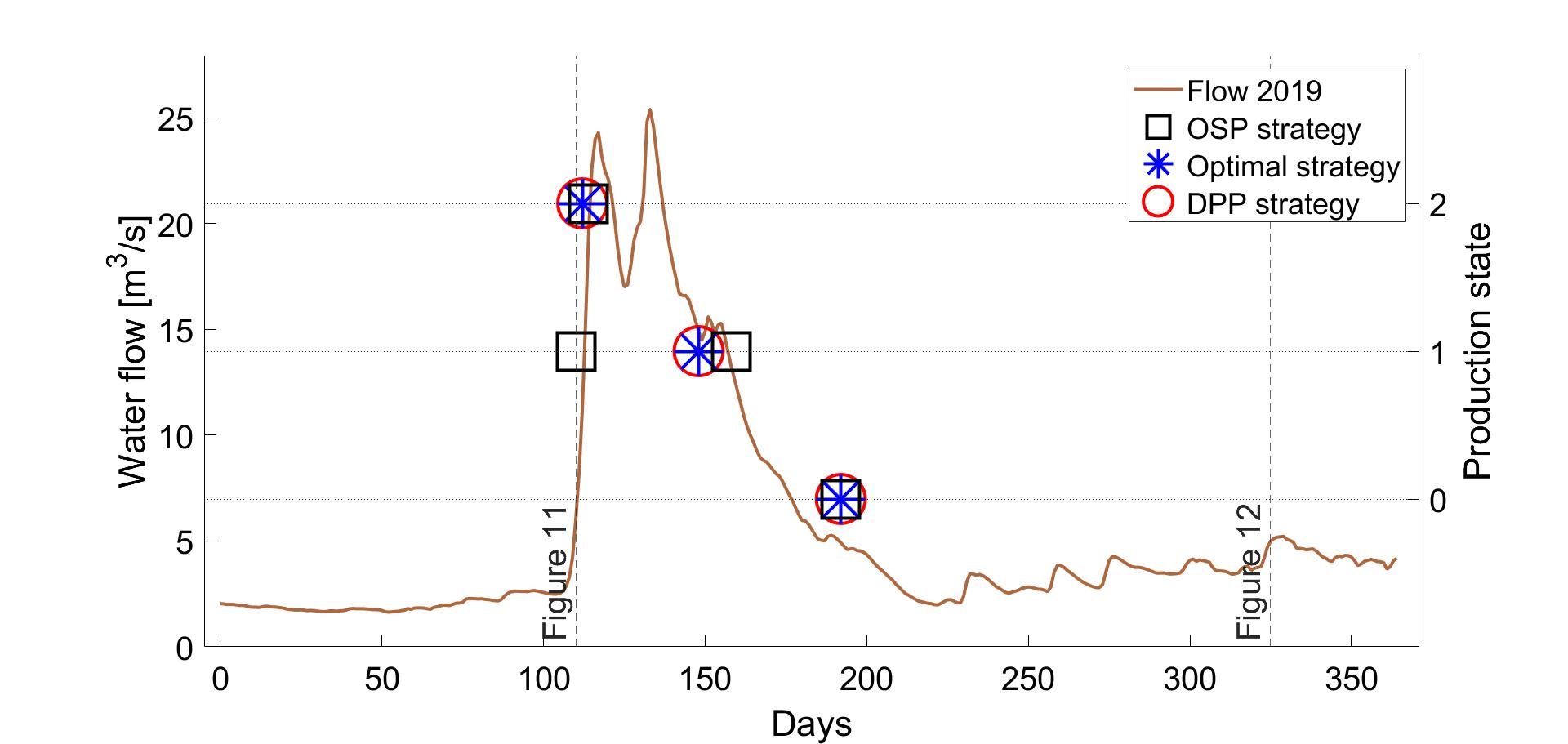}
    \caption{2019 with 5 days forecast and $C/D=0.0075$. There are discrepancies between the strategies suggested by OSP and DPP; the OSP takes a ''safer route" whereas the DPP strategy is more aggressive and finds the true optimum.} 
    \label{fig:stratdiff2019}
\end{figure}

  \begin{figure}
    \centering
    \includegraphics[width=\textwidth]{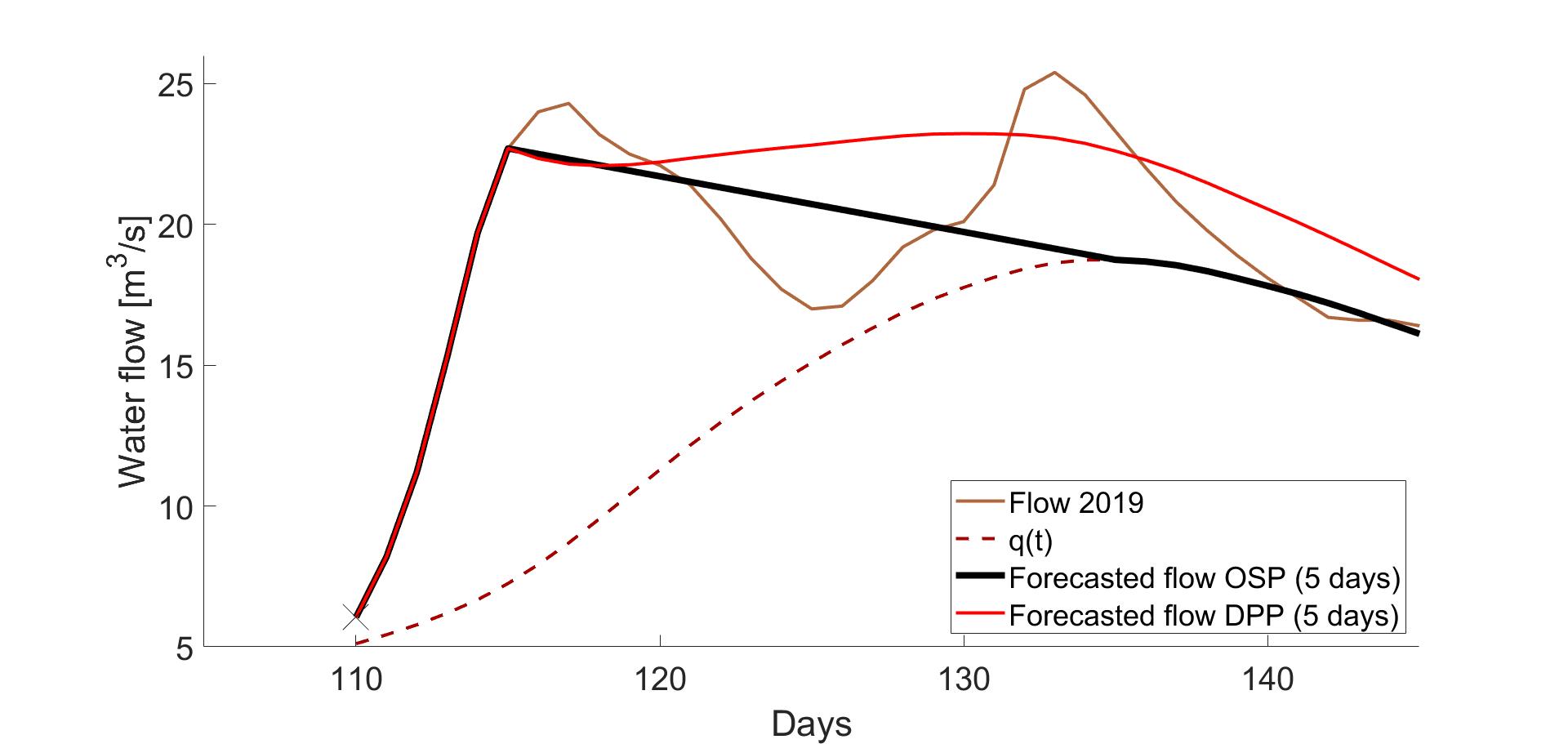}
    \caption{The OSP forecasts a lower flow than the DPP, therefore sub-optimally starting production in the intermediate mode $i=1$ at $t=110$.}
    \label{fig:closerinspection2019a}
\end{figure}

    \begin{figure}
    \centering
\includegraphics[width=\textwidth]{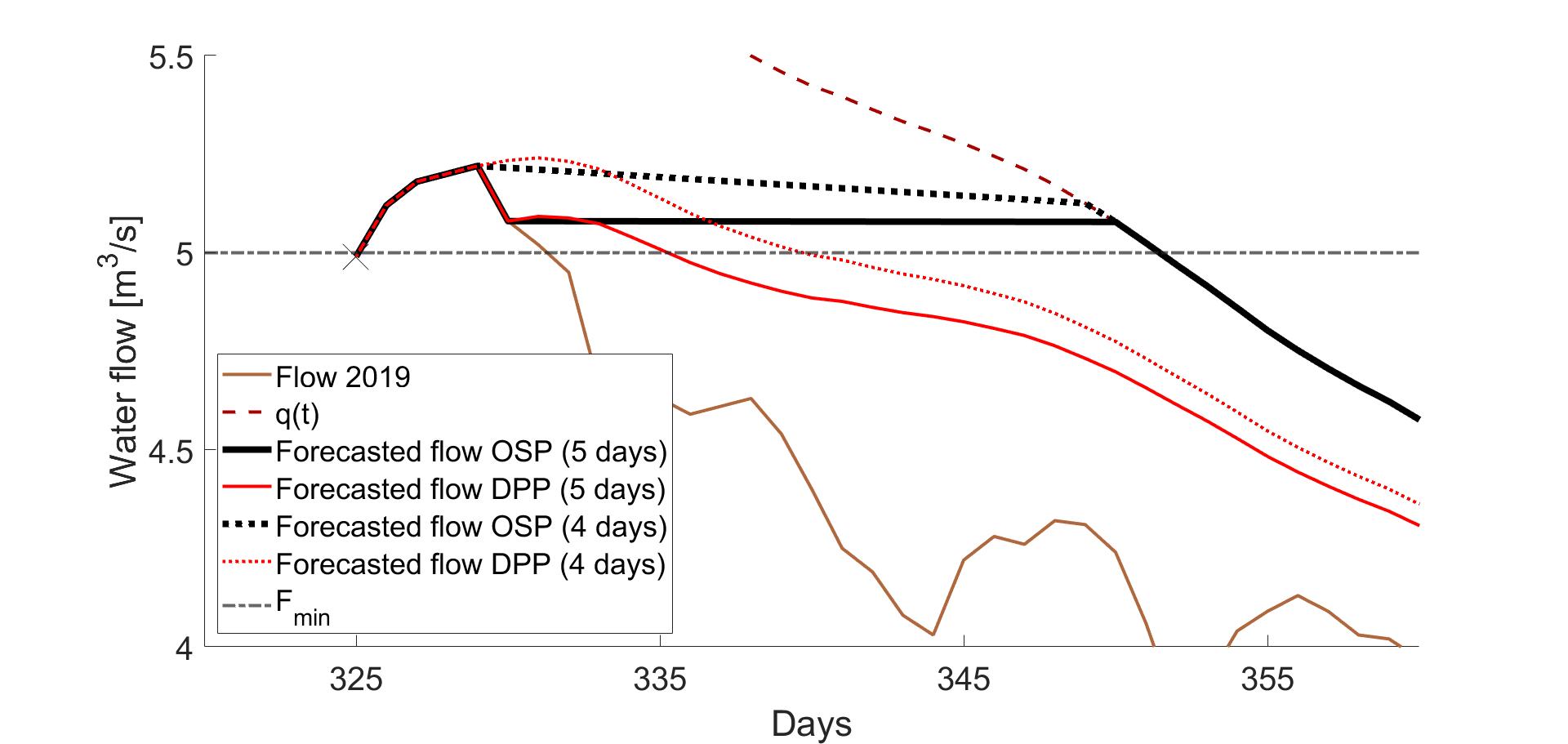}
    \caption{Both methods ignore opening production at $t=325$ with a $5$ day forecast, despite the flow being over the critical threshold $F_{min}=5 \,m^3/s$. With a shorter forecast of $4$ days, the projected flow is slightly higher and the DPP sub-optimally starts production while the OSP method still makes the correct decision and avoids costly opening and closing of the plant.}
    \label{fig:closerinspection2019b}
\end{figure}

\begin{remark}
\label{remark:OSPflow}
Note that the flow model of \cite{OOL22}, building on stochastic differential equations, is more complicated than the deterministic approach used in this paper. However, it is a simple task to adapt that flow model to meet the requirements of the deterministic approach presented in the current paper: simply set $\sigma=0$ in equation (3.1) of \cite{OOL22}. 
This gives a deterministic flow model that can be used 
as outlined above. When combining that flow model with the deterministic DPP-approach suggested in this paper some minor differences can be observed in the results, but the general qualitative observations made in Section \ref{sec:oldsetup} remain valid. This indicates that the simple flow model \eqref{eq:flowmodel_new} presented in Section \ref{sec:river-model} is sufficiently rich to tackle the problem at hand. This stays true also when decreasing or increasing the half-life $T_{1/2}$ in \eqref{eq:flowmodel_new}. 
In particular, performing the calculations of Section \ref{sec:results} for $T_{1/2} = 5$ and $T_{1/2} = 20$ gives an average of $97.2$ \% and $97.7$ \% of the theoretical maximum, respectively, for the years 2015-2022.
     \end{remark}

\section{Discussion}
\label{sec:discussionanconclusions}

The major upside of the scheme presented here is its simplicity, both from a mathematical and modeling perspective. The method can easily be adjusted and expanded to more complicated power plants with a large number of different production modes. In comparison, expanding the optimal switching-based model of \cite{OOL22} to include dams as in Section \ref{sec:method} would require treatment of interconnected PDEs with Neumann boundary conditions and computational proficiency to solve these explicitly. Moreover, the addition of further underlying processes in that setting increases the dimensionality of the underlying PDE and quickly requires explicit solutions of high-dimensional PDEs to be computationally tractable.
On the contrary, the computational resources needed here are relatively small and the method can cope with a larger number of underlying processes, at least as long as these are truly exogenous, e.g., wind, water flow, electricity price, etc. Moreover, the number of production modes can be increased without slowing down the process notably.
The number of processes that are affected by our actions must however be in the low single digits for the method to be tractable as we must keep track of all possible choices for these in the backward recursion. Relying on a coarse discretization and interpolation could push this limit a bit but at the risk of losing accuracy.

The downside of the method is primarily its lack of stochastic features, meaning that all uncertainty must be considered in the respective models for the underlying processes, possibly aggravating the modeling at that stage.  Moreover, a deterministic approach can lead to too many actions when the underlying process fluctuates around key values (e.g., $F_{min}$ in the example of Section \ref{sec:oldsetup}). This does not seem to be the case in our examples, but it must nevertheless be considered and observed closely in applications. 

The presented method uses deterministic dynamic programming but relates to stochastic dynamic programming for reservoir operations \cite{SDP1}. The key difference is that, rather than representing the possible states with a probability distribution, we construct and use deterministic models (see Section \ref{sec:river-model}) for characterization of the future states by simplifying historical data to its seasonal mean. 
In particular, the half-life reversion toward seasonal flow is a very simple model that efficiently merges short-term forecasts with historical data. 
As observed already in Remark~\ref{remark:OSPflow}, our results are not sensitive to the choice of $T_{1/2}$, indicating that this simple idea is already enough -- it is most likely not much more to be gained from polishing the model as long as the ingoing parameters are \textit{reasonable}. What is important in the model seems to be the combination of short-term forecast and historical mean flow, not the (often difficult) exact calibration of the parameters.

When it comes to comparison of the results of the DPP and the optimal switching-based method, the OSP is more stable and performs better with less accurate data,  but rarely finds the true optimal strategy. This is due to a conceptual difference between the DPP method suggested here and that based on stochastic differential equations in \cite{OOL22}; the former takes the input flow model as a 'fait accompli' while the latter expects stochastic deviations and tries to maximize the \textit{expected} profit. The OSP schemes therefore ''wait and see", thereby missing the true optimum slightly at the benefit of minimizing the risk of costly switches back and forth. With reliable forecasts and/or low costs of switching it therefore seems reasonable to opt for deterministic methods such as that suggested here while stochastic features are advantageous when information is scarce or insecure. 

Hydropower is often presented as a clean and renewable energy source that is environmentally preferable to fossil fuels or nuclear power. However, it often transforms rivers 
by, e.g., 
reducing flow velocity 
and disrupting sediment dynamics, and 
by extension, it therefore also alters riverine biodiversity. Freshwater ecosystems are in fact among the world's most threatened ecosystems \cite{sweden, YT22}. Therefore, an important challenge for river management is to identify situations where measures involving relatively small production losses can have major ecological advantages.
This calls for an extension of the present work towards a multi-objective optimization approach in which one imposes restrictions on, e.g., the reservoir level and the output flow from the power plant. A suggested strategy would in that case not consist of a single action but rather a Pareto-front consisting of efficient strategies where the manager can make a choice depending on the desired degree of environmental friendliness. In such multi-objective optimization, the simplicity of the current scheme can be a great advantage as it eases the addition of further traits for consideration.

\section*{Acknowledgments}
We would like to thank four anonymous reviewers for valuable comments and suggestions that has improved this manuscript. André Berg gratefully acknowledge the financial support provided by the Kempe Foundation.

\section*{Data availability statement}
The hydrological data used in this study is openly available via the Swedish Meteorological and Hydrological Institute at \url{https://vattenwebb.smhi.se/hydronu/} (Reference 24099, ''Sävarån").









\medskip
Received xxxx 20xx; revised xxxx 20xx; early access xxxx 20xx.
\medskip


\begin{thebibliography}{99}


\bibitem{DynProg86}
Allen R. B., Bridgeman S. G., 
{\it Dynamic programming in hydropower scheduling.}
Journal of Water Resources Planning and Management 112.3 (1986): 339--353.



\bibitem{bok}
Basson M.S., Allen R.B., Pegram G.G.S., van Rooyen  J.A.,
{\it Probabilistic management of water resource and hydropower systems.} Water Resources Publications, 1994



\bibitem{C_etal15}
Cheng C., et al., 
{\it China's small hydropower and its dispatching management.} 
Renewable and Sustainable Energy Reviews 42 (2015): 43--55.


\bibitem{DFHD21}
Danso D. K., François B., Hingray B.,  Diedhiou A.  
{\it Assessing hydropower flexibility for integrating solar and wind energy in West Africa using dynamic programming and sensitivity analysis. Illustration with the Akosombo reservoir, Ghana.} 
Journal of Cleaner Production, 287  (2021): 125559.



\bibitem{DWP19}
Dobson B., Wagener T., Pianosi F.,
{\it An argument-driven classification and comparison of reservoir operation optimization methods.} 
Advances in Water Resources, 128 (2019): 74--86.




\bibitem{IDynProg17}
Feng Z.-k., 
Niu W.-j.,
Cheng C.-t.
Wu X.-y.,
{\it Optimization of large-scale hydropower system peak operation with hybrid dynamic programming and domain knowledge.} 
Journal of Cleaner Production 171 (2018): 390--402.


\bibitem{rew21}
Giuliani M., 
Lamontagne J. R., 
Reed P. M., 
Castelletti A. 
{\it A State‐of‐the‐art review of optimal reservoir control for managing conflicting demands in a changing world.} 
Water Resources Research, 57(12) (2021), e2021WR029927.

\bibitem{SDP2}
Gjelsvik A., Mo B.,  Haugstad A., 
{\it Long-and medium-term operations planning and stochastic modelling in hydro-dominated power systems based on stochastic dual dynamic programming.} 
Handbook of power systems I (2010): 33-55.

\bibitem{rew04}
Labadie J. W. 
{\it Optimal operation of multireservoir systems: State-of-the-art review.} 
Journal of water resources planning and management, 130(2), (2004), 93--111.


\bibitem{highdim}
Mitjana F., Denault M.,  Demeester, K. 
{\it Managing chance-constrained hydropower with reinforcement learning and backoffs.} 
Advances in Water Resources, 169  (2022): 104308.

\bibitem{NS97}
Nilsson O.,  
Sjelvgren D. 
{\it Hydro unit start-up costs and their impact on the short term scheduling strategies of Swedish power producers.} 
IEEE Transactions on power systems, 12(1), (1997), 38--44.


\bibitem{OOL22}
Olofsson M.,  \"Onskog T., Lundstr\"om N.L.P., 
{\it Management strategies for run-of-river hydropower plants: an optimal switching approach.} 
Optimization and Engineering (2021): 1--25.

\bibitem{Queiroz}
de Queiroz A. R. {\it Stochastic hydro-thermal scheduling optimization: An overview.} 
Renewable and Sustainable Energy Reviews 62 (2016): 382-395.




\bibitem{sweden}
Ren\"of\"alt Malm B., Jansson R., Nilsson C.,
{\it Effects of hydropower generation and opportunities for environmental flow management in Swedish riverine ecosystems.}
Freshwater Biology 55.1 (2010): 49--67.


\bibitem{singh}
Singh V. K., Singal S. K., 
{\it Operation of hydro power plants--a review.} Renewable and Sustainable Energy Reviews 69 (2017): 610-619.



\bibitem{pump}
Wamalwa F., Sichilalu S.,  Xia X., 
{\it Optimal control of conventional hydropower plant retrofitted with a cascaded pumpback system powered by an on-site hydrokinetic system.} 
Energy Conversion and Management 132 (2017): 438--451.



\bibitem{Y10}
Yang J. 
{\it Underwater tunnel piercing in refurbishment of Akkats power station.} 
ICOLD Annual Meeting and Symposium, May (2010), Hanoi, Vietnam


\bibitem{Y18}
Yang J., Andreasson P., H\"ogstr\"m C.-M., Teng P. 
{\it The tale of an intake vortex and its mitigation
countermeasure: a case study from Akkats hydropower station.} 
Water 10(7) (2018): 881



\bibitem{YT22}
Yoshioka H., Tsujimura M., 
{\it Hamilton-Jacobi-Bellman-Isaacs equation for rational inattention in the long-run management of river environments under uncertainty}. 
Computers and Mathematics with Applications, 112 (2022): 23--54.




\bibitem{YSA05}
Yurtal R., Seckin G., Ardiclioglu G. 
{\it Hydropower optimization for the lower Seyhan system in Turkey using dynamic programming.} 
Water international, 30(4) (2005): 522--529.


\bibitem{SDP1}
Zambelli M. S., S. Soares, D. da Silva. 
{\it Deterministic versus stochastic dynamic programming for long term hydropower scheduling.} 
2011 IEEE Trondheim PowerTech. IEEE, 2011.




\bibitem{IDynProg14}
Zhao T., Zhao J., Yang D., 
{\it Improved dynamic programming for hydropower reservoir operation.} Journal of Water Resources Planning and Management 140.3 (2014): 365--374.

\end{thebibliography}
\end{document}